\nonstopmode
\documentclass[10pt,twocolumn,smallcolumn]{quick-article}

\renewcommand{\ShowTODO}[1]{}

\usepackage{fancyhdr}

\makeatletter
\font\tensmc=cmcsc10

\def\footnotesize{\@setsize\footnotesize{10pt}\viiipt\@viiipt
        \abovedisplayskip \z@
        \belowdisplayskip\z@
        \abovedisplayshortskip\abovedisplayskip
        \belowdisplayshortskip\belowdisplayshortskip
  \def\@listi{\leftmargin\leftmargini \topsep 3pt plus 1pt minus 1pt
        \parsep 2pt plus 1pt minus 1pt
        \itemsep \parsep}}

\skip\footins 12pt plus 12pt

\def\footnoterule{\kern3\p@  \hrule width 3em} % the \hrule is .4pt high

\def\ps@plain{\let\@mkboth\@gobbletwo
     \def\@oddfoot{{\hfil\small\thepage\hfil}}%
     \def\@oddhead{}
      \def\@evenhead{}\def\@evenfoot{}}

\def\ps@headings{\let\@mkboth\markboth
        \def\@oddfoot{}\def\@evenfoot{}%
        \def\@evenhead{{\rm\thepage}\hfil{\small\leftmark}}%
        \def\@oddhead{{\noindent\small\rightmark}\hfil{\rm\thepage}}%

\def\ps@myheadings{\let\@mkboth\@gobbletwo
 \def\@oddfoot{}\def\@evenfoot{}%
 \def\@oddhead{\rlap{\normalsize\rm\rightmark}\hfil{small\thepage}}%
 \def\@evenhead%{\hfil{\small\@chapapp}\
                {\small\thepage}\hfil\llap{\normalsize\rm\leftmark}}%
        \def\chaptermark##1{}%
       \def\sectionmark##1{}\def\subsectionmark##1{}}

\if@titlepage
  \renewcommand\maketitle{\begin{titlepage}%
  \let\footnotesize\small
  \let\footnoterule\relax
  \let \footnote \thanks
  \null\vfil
  \vskip 60\p@
  \begin{center}%
    {\LARGE \@title \par}%
    \vskip 3em%
    {\large
     \lineskip .75em%
      \begin{tabular}[t]{c}%
        \@author
      \end{tabular}\par}%
      \vskip 1.5em%
    {\large \@date \par}%       % Set date in \large size.
  \end{center}\par
  \@thanks
  \vfil\null
  \end{titlepage}%
  \setcounter{footnote}{0}%
  \global\let\thanks\relax
  \global\let\maketitle\relax
  \global\let\@thanks\@empty
  \global\let\@author\@empty
  \global\let\@date\@empty
  \global\let\@title\@empty
  \global\let\title\relax
  \global\let\author\relax
  \global\let\date\relax
  \global\let\and\relax
}
\else
\renewcommand\maketitle{\par
  \begingroup
    \renewcommand\thefootnote{\@fnsymbol\c@footnote}%
    \def\@makefnmark{\rlap{\@textsuperscript{\normalfont\@thefnmark}}}%
    \long\def\@makefntext##1{\parindent 1em\noindent
            \hb@xt@1.8em{%
                \hss\@textsuperscript{\normalfont\@thefnmark}}##1}%
    \if@twocolumn
      \ifnum \col@number=\@ne
        \@maketitle
      \else
        \twocolumn[\@maketitle]%
      \fi
    \else
      \newpage
      \global\@topnum\z@   % Prevents figures from going at top of page.
      \@maketitle
    \fi
    \thispagestyle{fancy}\@thanks
  \endgroup
  \setcounter{footnote}{0}%
  \global\let\thanks\relax
  \global\let\maketitle\relax
  \global\let\@maketitle\relax
  \global\let\@thanks\@empty
  \global\let\@author\@empty
  \global\let\@date\@empty
  \global\let\@title\@empty
  \global\let\title\relax
  \global\let\author\relax
  \global\let\date\relax
  \global\let\and\relax
}
\def\@maketitle{%
  \newpage
  \null
  \vskip 2em%
  \begin{center}%
  \let \footnote \thanks
    {\LARGE \@title \par}%
    \vskip 1.5em%
    {\large
      \lineskip .5em%
      \begin{tabular}[t]{c}%
        \@author
      \end{tabular}\par}%
    \vskip 1em%
    {\large \@date}%
  \end{center}%
  \par
  \vskip 1.5em}
\fi

\long\def\appendix{\par\setcounter{section}{0}\setcounter{subsection}{0}\setcounter{equation}{0}\gdef\theequation{\@Alph\c@section.\arabic{equation}}%
\gdef\thesection{\@Alph \c@section }}

\makeatother

\title{
    Split-Decomposition Trees with Prime Nodes:\\
    Enumeration and Random Generation of Cactus Graphs}
\author{%
  Maryam Bahrani\thanks{Department~of~Computer~Science, Princeton
    University, 35~Olden~Street, Princeton, NJ 08540, USA.
    \email{mbahrani@princeton.edu} and \email{lumbroso@cs.princeton.edu}}\and %
  J\'{e}r\'{e}mie Lumbroso\footnotemark[1]}
\date{}

\newcommand{\TODOAll}[1]{\TODOTemplate{Orchid}{{\bf TODO:} #1}}%
\newcommand{\TODOMaryam}[1]{\TODOTemplate{Periwinkle}{{\bf TODO Maryam:} #1}}%
\newcommand{\TODOJeremie}[1]{\TODOTemplate{LimeGreen}{{\bf TODO J\'er\'emie:} #1}}%

\def\mLeaf{\bullet}    % a leaf
\combstruct{Useq}
\combstruct{Ucyc}
\combstruct{Xxx}
\combstruct{Yyy}
\combstruct{Zzz}
\combstruct{Www}

\def\mLeaf{\bullet}    % a leaf

\def\clsFRG{\cls[\mathrm{fr}]{G}}
\def\clsPRG{\cls[\mathrm{pr}]{G}}
\def\clsFUG{\cls[\mathrm{fu}]{G}}
\def\clsPUG{\cls[\mathrm{pu}]{G}}

\begin{document}
\maketitle

\fancyfoot[R]{\footnotesize{\textbf{Copyright \textcopyright\ 2018 by SIAM\\
Unauthorized reproduction of this article is prohibited}}}

\section*{Abstract}

In this paper, we build on recent results by Chauve \textit{et~al.} and Bahrani and Lumbroso, which combined the split-decomposition, as exposed by Gioan and Paul, with analytic combinatorics, to produce new enumerative results on graphs---in particular the enumeration of several subclasses of perfect graphs (distance-hereditary, 3-leaf power, ptolemaic).

Our goal was to study a simple family of graphs, of which the split-decomposition trees have prime nodes drawn from an enumerable (and manageable!) set of graphs. Cactus graphs, which we describe in more detail further down in this paper, can be thought of as trees with their edges replaced by cycles (of arbitrary lengths). Their split-decomposition trees contain prime nodes that are cycles, making them ideal to study.

We derive a characterization for the split-decomposition trees of cactus graphs, produce a general template of symbolic grammars for cactus graphs, and implement random generation for these graphs, building on work by Iriza.

%\end{rewrite}

\begin{figure*}[p!]
  \centering
  \begin{minipage}[t]{.65\linewidth}
    \centering
    \includegraphics[scale=0.8, page=2]{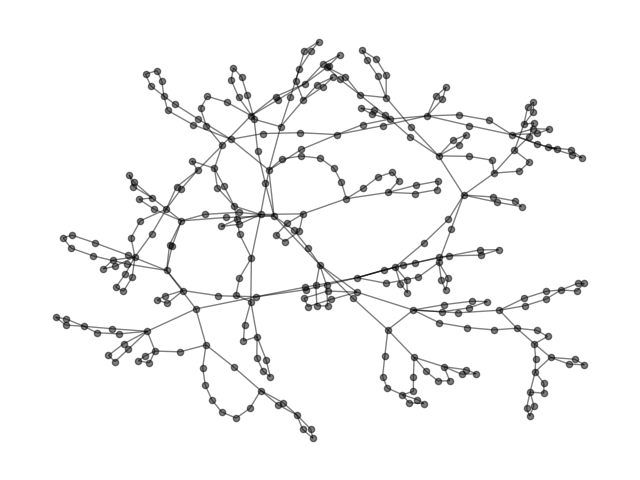}
    \subcaption{\label{fig:random-cactus-small}%
      A random mixed cactus graph with 309 vertices and 80 cycles.}
    \hspace{6pc}%
    \includegraphics[scale=0.8, page=3]{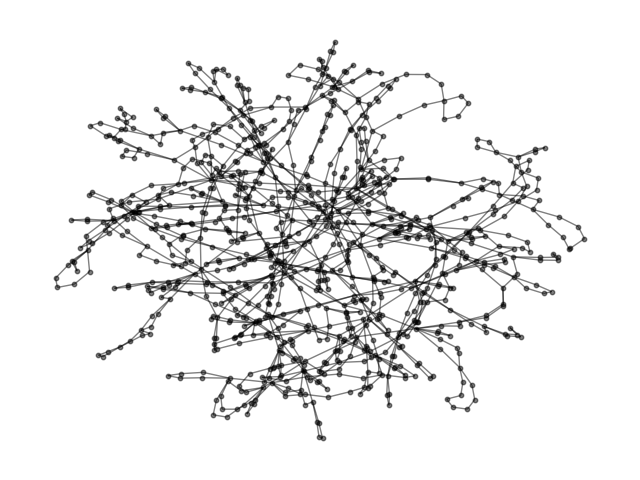}
    \subcaption{\label{fig:random-cactus-large}%
      A random mixed cactus graph with 933 vertices and 239 cycles.}
  \end{minipage}\hspace{0.08\linewidth}%
  \caption{\label{fig:random-cacti} %
    Our symbolic grammars can be used to create Boltzmann samplers for the uniform random generation of large cactus graphs. Here we have included two cactus graphs, each of which has been drawn uniformly at random from the family of unlabeled, plane, rooted cacti where each cycle has at least 4 nodes.}
\end{figure*}

\section*{Introduction}

\TODOJeremie{Write segue for split-decomposition using tree decompositions.}

In general, counting graphs is much more difficult than counting trees---the latter are completely recursive, and as such, there are many generic results and theorems that can be applied straightforwardly~\cite{FlSe09}. Graphs, on the other hand, can have cycles and can thus not be specified recursively in the general case. A \emph{tree decomposition} is a method by which graphs are put in correspondence with rooted trees (usually bijectively), which are easier to specify and study; to convert a graph to a tree, the decomposition makes certain hypotheses on the structure of the graphs.

There exists a number of tree decompositions, and several have been used in the context of analytic combinatorics. Most recently, Chauve~\textit{et~al.}~\cite{ChFuLu17} have used the split decomposition to provide the first exact enumeration of distance-hereditary graphs (the largest subset of graphs that can be entirely decomposed with the split-decomposition).  
\subsection*{The split decomposition and prime nodes.}
%\begin{rewrite}
Informally, a \emph{split} in a graph is a bipartition of the vertices into two sets of size at least 2, such that the edges between the sets form a complete bipartite graph.  All splits of a graph form a tree-like structure, where each internal node of the tree is labeled with a graph. The split decomposition~\cite{GiPa12} refers to the process of transforming a graph into a \emph{graph-labeled tree} using splits of the original graph.

The split decomposition was first introduced by Cunningham~\cite{Cunningham82} and later reformulated by Gioan and Paul~\cite{GiPa12}, who introduced the concept of \emph{graph-labeled trees}. More recently, Chauve~\textit{et al.}~\cite{ChFuLu14, ChFuLu17} observed that these graph-labeled trees are well-suited for recursive symbolic specification, which they used to provide, among other things, the first exact enumeration of distance-hereditary graphs. Distance-hereditary graphs are the totally decomposable graphs for the split decomposition~\cite{GiPa12,Cunningham82}, and therefore a natural first class of graphs to study when considering the split decomposition for enumerating graphs.

The paper by Chauve~\etal~\cite{ChFuLu17} was first in a line of work on the application of the split decomposition to graph enumeration. Bahrani and Lumbroso~\cite{BaLu17} used the same framework to give grammars for \emph{subsets} of distance-hereditary graphs that are specified in terms of forbidden induced subgraphs. They showed that constraining a split decomposition tree to avoid certain patterns can help avoid corresponding induced subgraphs in the original graph. The ability to \emph{forbid} induced subgraphs led to grammars and full enumeration for a wide set of graph classes: ptolemaic, block, and variants of cactus graphs (2,3-cacti, 3-cacti and 4-cacti). In certain cases, no enumeration was known (ptolemaic, 4-cacti); in other cases, although the enumerations were known, an abundant potential is unlocked by the provided grammars (in terms of asymptotic analysis, random generation, and
parameter analysis, \textit{etc.}).

The methodology has shown its potential in this line of work, specially in providing grammars for totally decomposable graph classes, \textit{i.e.} distance-hereditary graphs and their subfamilies. As stated by Bahrani and Lumbroso~\cite{BaLu17}, the natural next step to develop its reach is to study split decomposition trees that are not fully decomposable and contain prime nodes. This is part of the motivation behind studying cactus graphs, which result in split decomposition trees with prime polygonal nodes. As such, this paper is the first work in the line of papers introduced above to examine split decomposition trees that are not fully decomposable.

\subsection*{Cactus graphs.}

\begin{figure*}[ht!]
    \centering
    \includegraphics[scale=0.3]{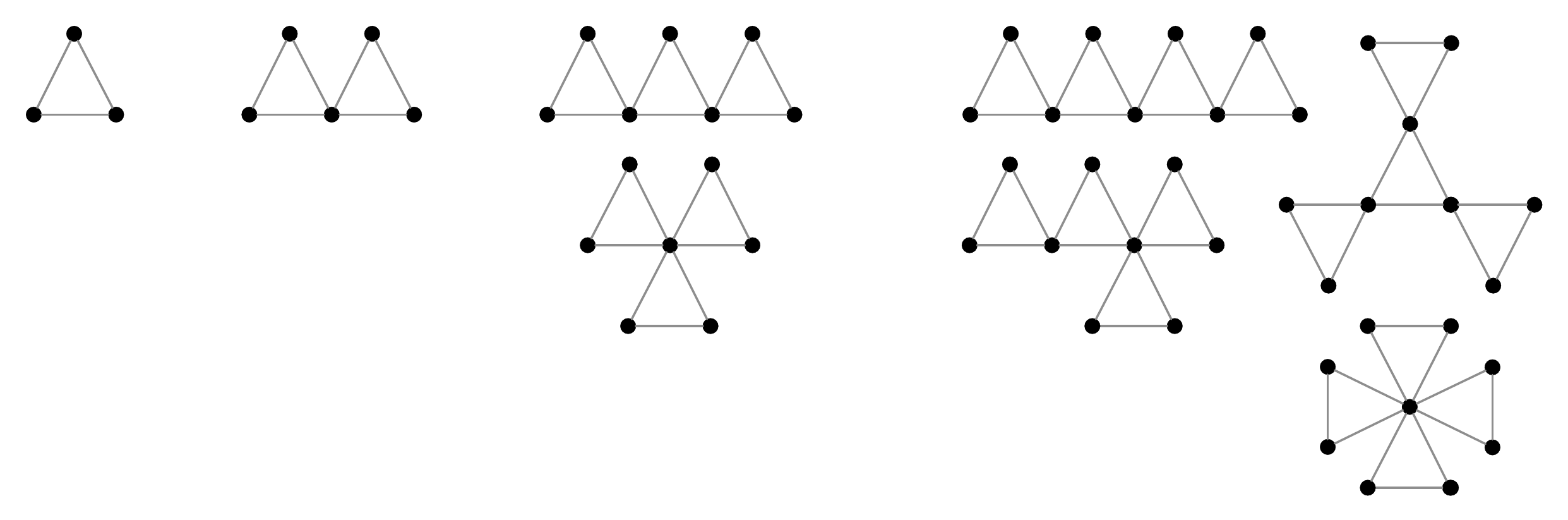}
    \caption{\label{fig:example-3cactus}%
      Small (unrooted, unlabeled, pure, free) 3-cacti}
\end{figure*}

A cactus is a graph in which no two cycles share an edge\footnote{In this paper, it is assumed that cactus graphs are connected.}. It can be thought of as a collection of polygons and bridges that pairwise share at most one vertex. Cactus graphs come in many flavors (\textit{e.g.} see Figure~\ref{fig:example-3cactus}):
\begin{itemize}[nosep, noitemsep]
    \item A cactus graph is \emph{pure} if all its cycles have the
    same size; it is \emph{mixed} otherwise.
    \item A graph is \emph{labeled} if its vertices are distinguished for the purpose of enumeration; it is \emph{unlabeled} if the vertices have no distinct identification besides their adjacencies.
    \item A \emph{rooted} graph has one distinguished vertex; it is otherwise called \emph{unrooted}.
    \item A \emph{plane} cactus graph is embedded in the plane so that every vertex is paired with a circular permutation of the edges incident to it; when this permutation is insignificant, the cactus graph is called \emph{non-plane} or \emph{free}.
\end{itemize}
\noindent
Besides being objects of interest in combinatorics, cactus graphs have applications in other fields. For example, many hard problems on general graphs have efficient algorithms on cactus graphs~\cite{baste2017parameterized,chandran2017algorithms,das2014cactus}. Cactus graphs are also used in modeling electronic circuits~\cite{nishi1984topological}, networks~\cite{arcak2011diagonal}, and comparative genomics~\cite{paten2011cactus}.

The problem of enumerating cactus graphs was first proposed in 1950 in a lecture by Uhlenbeck~\cite{UhlenbeckLecture}, where these graphs were referred to as \emph{Husimi Trees}. Subsequently, Harry and Uhlenbeck~\cite{HaUh53} as well as Ford and coauthors~\cite{ford1956combinatorialI,ford1956combinatorialII,ford1956combinatorialIII,ford1957combinatorialIV} and Bergeron~\textit{et~al.}~\cite{bergeron1998combinatorial} have examined the enumeration of unlabeled and labeled cactus graphs in a series of notes. These works derive functional equations for the generating functions of cactus graphs, with no convenient way of obtaining a few hundred terms of the enumeration or the prospect of building random generators. Additionally, their methodology is dependent upon the symmetric structure of cycles in cactus graphs and therefore hard to generalize to other graph classes. More recently, Bona \textit{et al.}~\cite{bona2000enumeration} derived the enumeration of \textit{plane} cactus graphs, but their methodology does not seem to generalize to non-plane cacti. 

The goal of this paper is to develop a simple, extendable framework that not only allows for enumerating all aforementioned varieties of cactus graphs but also is amenable to imposing arbitrary restrictions on the structure of the graphs\footnote{For example, one could enumerate cactus graphs where the length of each cycle is a prime number.} as well as random generation. This goal is achieved by developing symbolic grammars for different varieties of cactus graphs.

Since cactus graphs are very similar to trees, the same grammars can be developed more directly by exploiting this apparent recursivity. However, while the entire machinery introduced in this paper is not necessary to obtain grammars for cactus graphs, these tools belong to a more general and powerful framework developed by Chauve~\etal~\cite{ChFuLu17} and extended by Bahrani and Lumbroso~\cite{BaLu17}. This framework is highly flexible and can be generalized to obtain grammars for many other classes of cactus graphs. We therefore believe it is informative to state our results within this framework, specially since our original investigation was to look at identifying characterizable families of graphs of which the split decomposition trees contain a manageable family of prime nodes (in this case: various cycles and sequences).

%\end{rewrite}

%%%%%%%%%%%%%%%%%%%%%%%%%%%%%%%%%%%%%%%%%%%%%%%%
%%% SECTION on Definitions and Preliminaries %%%

\section{Definitions and preliminaries\label{sec:preliminaries}}

In this section, we introduce standard definitions from graph
theory (\ref{subs:def-graph}), followed by a formal introduction to the split-decomposition expressed in terms of graph-labeled trees  (\ref{subs:graph-labeled} and \ref{subs:split}), which are the tools also used by Chauve~\etal~\cite{ChFuLu17} and Bahrani and Lumbroso~\cite{BaLu17}.

\subsection{Graph definitions.\label{subs:def-graph}}

For a graph $G$, we denote by $V(G)$ its vertex set and $E(G)$ its edge set.

Moreover, for a vertex $x$ of a graph $G$, we denote by $N(x)$ the neighborhood of $x$, that is the set of vertices $y$ such that $\{x,y\}\in E(G)$; this notion extends naturally to vertex sets: if $\cramped{V_1}\subseteq V(G)$, then $N(\cramped{V_1})$ is the set of vertices defined by the (non-disjoint) union of the neighborhoods of the vertices in $\cramped{V_1}$. The subgraph of $G$ induced by a subset $\cramped{V_1}$ of vertices is denoted by $G[\cramped{V_1}]$.

A clique on $k$ vertices, denoted by $\cramped{K_k}$, is the complete graph on $k$ vertices (\textit{i.e.}, there exists an edge between every pair of vertices). A star on $k$ vertices, denoted by $\cramped{S_k}$, is the graph with one vertex of degree $k-1$ (the \emph{center} of the star) and $k-1$ vertices of degree $1$ (the \emph{extremities} of the star).

A closed walk in a graph is an alternating sequence of vertices and edges, starting and ending at the same vertex, where each edge is adjacent in the sequence to its endpoints. A cycle is a closed walk in which no repetitions of vertices and edges are allowed, except for the starting and ending point. We use the terms \emph{cycles} and \emph{polygons} interchangeably. In this paper we work with simple graphs (\textit{i.e.} graphs without self-loops or multi-edges), so all cycles have at least 3 distinct vertices.

\subsection{Graph-labeled trees.\label{subs:graph-labeled}}
We first introduce the notion of \emph{graph-labeled tree}, due to Gioan and Paul~\cite{GiPa12}, then define the split-decomposition and finally give the characterization of a \emph{reduced} split-decomposition tree, described as a graph-labeled tree.
\newpage
\begin{definition}\label{def:glt}
  A graph-labeled tree $(T,\cls{F})$ is a tree $T$ in which every internal node $v$ of degree $k$ is labeled by a graph $\cramped{G_v} \in \cls{F} $ on $k$ vertices, called \emph{marker vertices}, such that there is a bijection $\cramped{\rho_v}$ between the edges of $T$ incident to $v$ and the vertices of $\cramped{G_v}$.
\end{definition}

\begin{figure*}[h!]
  \centering
  \begin{bigcenter}
    \begin{minipage}[t]{.45\linewidth}
      \centering
      \includegraphics[scale=0.4]{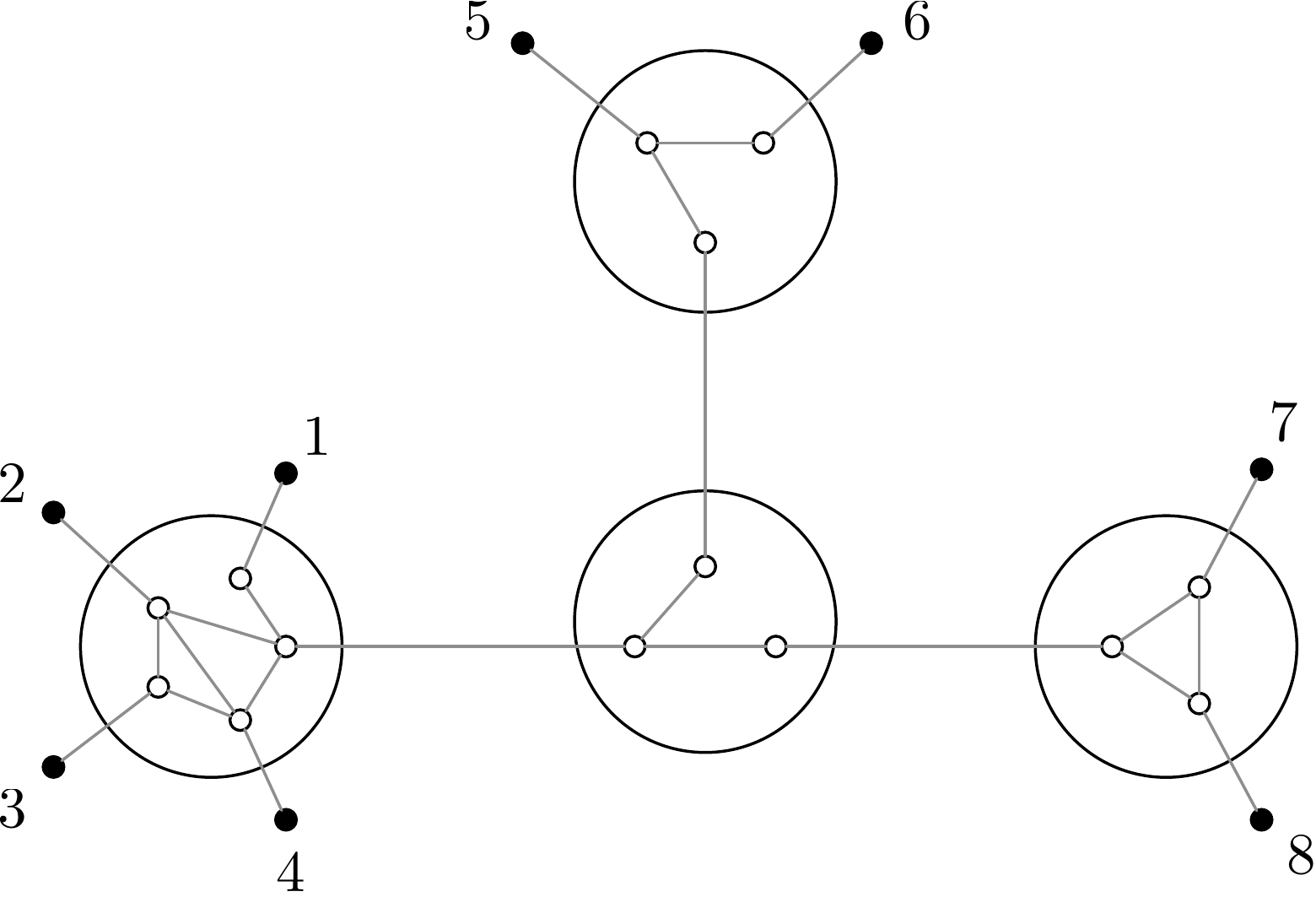}%
      \subcaption{\label{fig:ex-split-glt}%
        A graph-labeled tree.}
    \end{minipage}\hspace{2em}
    \begin{minipage}[t]{.45\linewidth}
      \centering
      \includegraphics[scale=0.4]{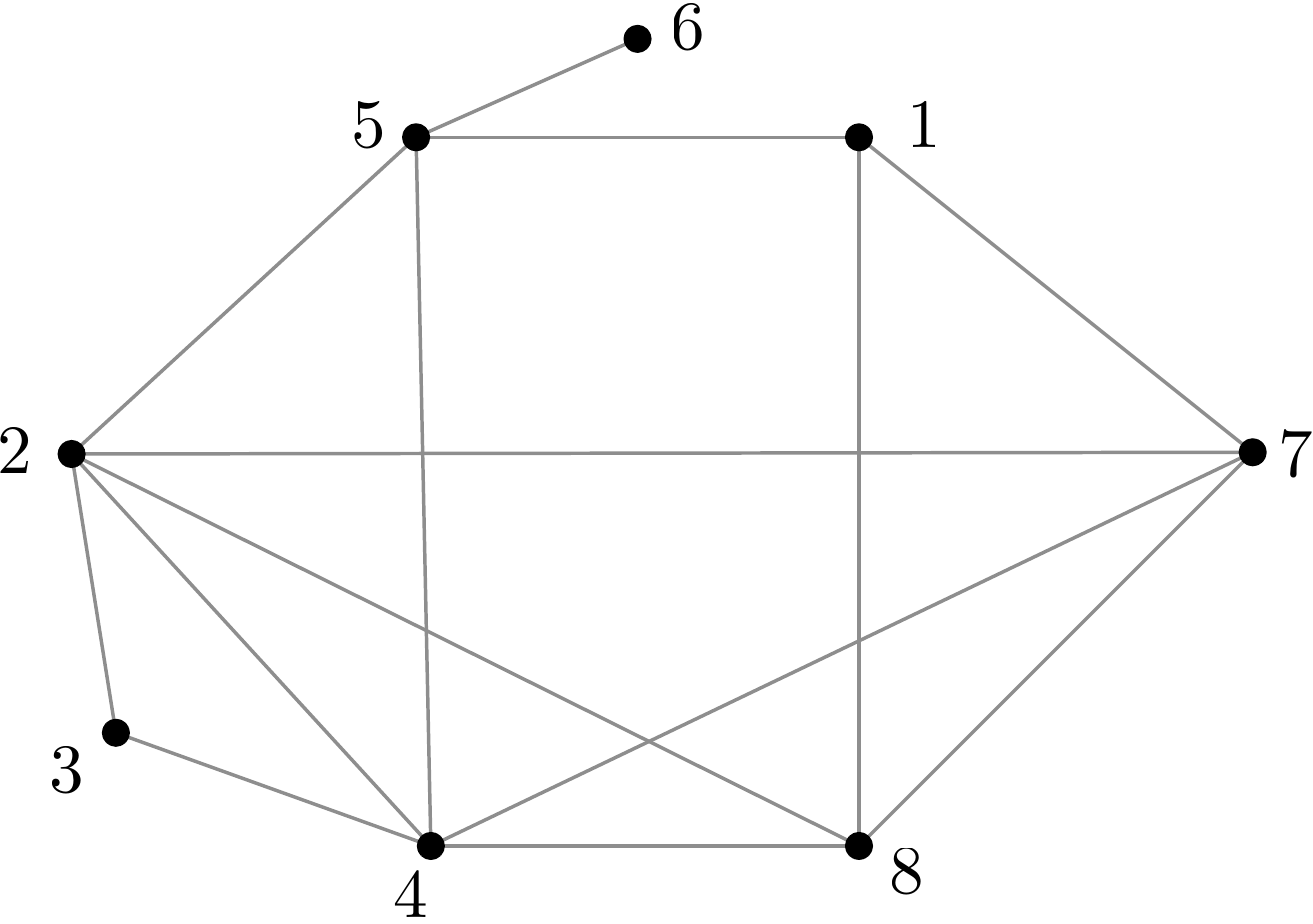}%
      \subcaption{\label{fig:ex-split-orig}%
        Original graph for (or \emph{accessibility graph} of)
        the graph-labeled tree in Figure~\ref{fig:ex-split-glt}.}
    \end{minipage}
  \end{bigcenter}
  \caption{%
    \label{fig:ex-split}%
    Two leaves of the split-decomposition graph-labeled tree (left)
    correspond to adjacent vertices in the original graph that was
    decomposed (right) if there exists an \emph{alternated path}: a path between those leaves, which uses at most one interior edge of any graph label. For example, vertex 5 is adjacent to vertex 4 in the original graph, because there is an alternated path between the two corresponding leaves in the split-decomposition tree; vertex 5 is not adjacent to vertex 3 however, because that would require the path to take two interior edges of the (\emph{prime}) leftmost graph-label.}
\end{figure*}

\noindent For example, in Figure~\ref{fig:ex-split} the internal nodes of $T$ are denoted by large circles, the marker vertices are denoted by small hollow circles, the leaves of $T$ are denoted by small solid circles, and the bijection $\cramped{\rho_v}$ is denoted by each edge that crosses the boundary of an internal node and ends at a marker vertex.

\begin{figure*}[h!]
  \centering
  \includegraphics[scale=0.55]{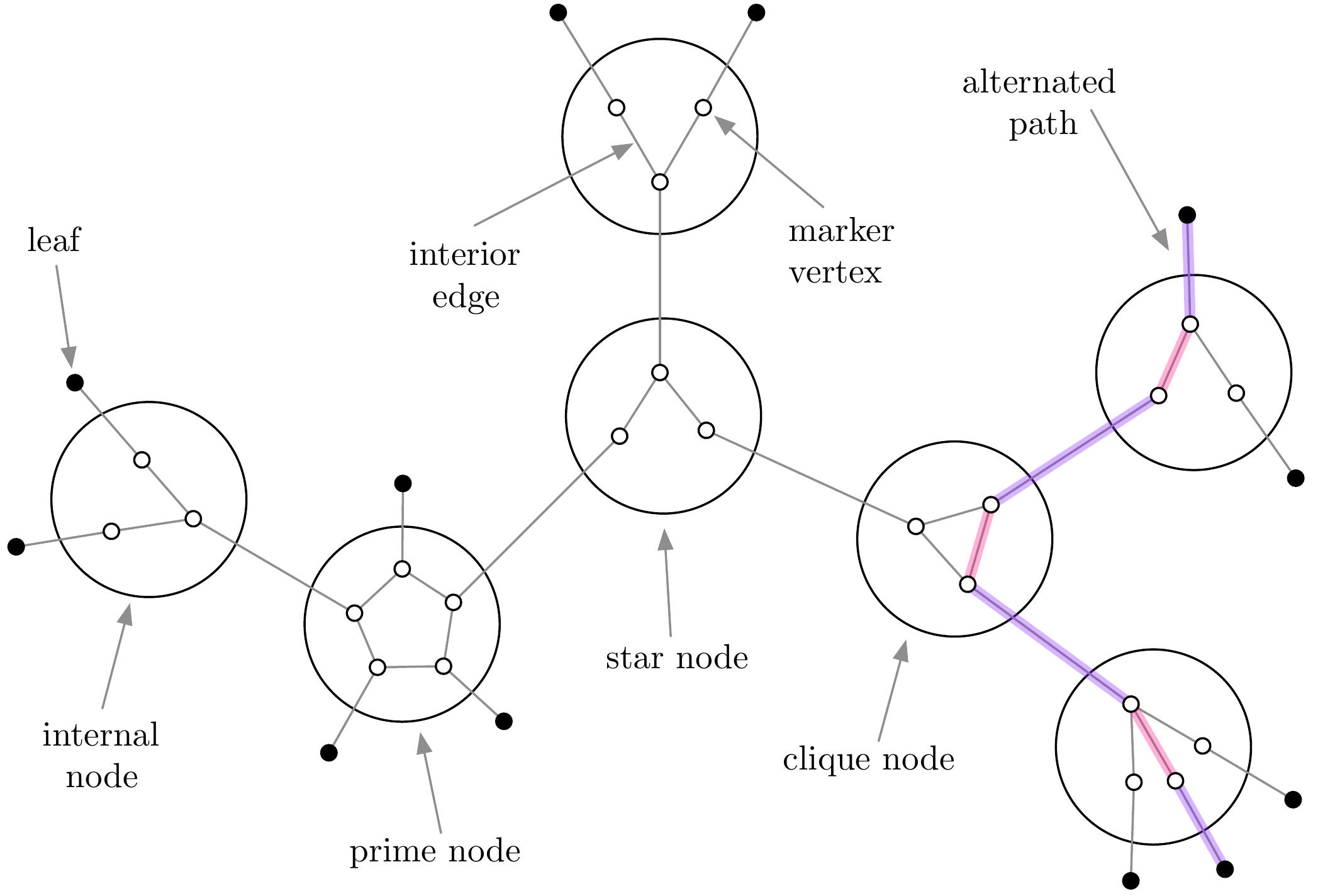}
  \caption{\label{fig:split-terminology}%
    In this figure, we present a few terms that we use a lot in this article.}
\end{figure*}

These graph-labeled trees are a powerful tool for studying the structure of the original graph they
describe. Some elements of the terminology have been summarized in
Figure~\ref{fig:split-terminology} (reproduced from Bahrani and Lumbroso~\cite{BaLu17}), as they are frequently referenced in the proofs of Section~\ref{sec:characterization}.

\begin{definition}
  Let $(T, \cls{F})$ be a graph-labeled tree and let $\ell, \ell'\in V(T)$ be leaves of $T$. We say that there is an \emph{alternated path} between $\ell$ and $\ell'$, if there exists a path from $\ell$ to $\ell'$ in $T$ such that for any adjacent edges $e = \set{u,v}$ and $e' = \set{v,w}$ on the path, $\set{\cramped{\rho_v}(e),\cramped{\rho_v}(e')}\in E(\cramped{G_v})$.
\end{definition}

\begin{definition}%
  \label{def:split-originalgraph}%
  The \emph{original graph}, also called \emph{accessibility graph}, of a graph-labeled tree $(T, \cls{F})$ is the graph $G$ where $V(G)$ is the leaf set of $T$ and, for $x, y\in V(G)$, $\set{x,y}\in E(G)$ iff there is an alternated path between $x$ and $y$ in $T$.
\end{definition}

\noindent Figures~\ref{fig:ex-split} and~\ref{fig:split-terminology} (reproduced from Bahrani and Lumbroso~\cite{BaLu17}) illustrate the concept of an alternated path: it is, more informally, a path that only ever uses at most one interior edge of any graph label.

%%%%%%%%%%%%%%%%%% FIGURES %%%%%%%%%%%%%%%%%%

%%%%%%%%%%%%%%%%%%%%%%%%%%%%%%%%%%%%%%%%%%%%%

\subsection{Split-decomposition.\label{subs:split}}
In this section we summarize various concepts, a more formal and comprehensive treatment of which can be found in previous work~\cite{ChFuLu17,BaLu17}.

A \emph{split} in a graph is a bipartition of the vertices into two sets of size at least 2, such that the edges between the two sets form a complete bipartite graph. More formally,
\begin{definition}%
  \label{def:split}%
  A \emph{split}~\cite{Cunningham82} of a graph $G$ with vertex set $V$ is a bipartition $(\cramped{V_1},\cramped{V_2})$ of $V$ (\textit{i.e.}, $V=\cramped{V_1}\cup V_2$, $\cramped{V_1}\cap \cramped{V_2}=\emptyset$) such that
  \begin{enumerate}[label=(\alph*), noitemsep, nosep]
  \item $|V_1|\geqslant 2$ and $|V_2|\geqslant 2$;
  \item every vertex of $N(\cramped{V_1})$ is adjacent to every vertex of $N(\cramped{V_2})$.
  \end{enumerate}
\end{definition}
\noindent
A graph without any split is called a \emph{prime} graph. A graph is \emph{degenerate} if any bipartition of its vertices into sets of size at least 2 is a split: cliques and stars are the only such graphs.

Informally, the split-decomposition of a graph $G$ consists of finding a split $(\cramped{V_1}, \cramped{V_2})$ in $G$, followed by decomposing $G$ into two graphs $\cramped{G_1}=G[\cramped{V_1}\cup \{\cramped{x_1}\}]$, where $\cramped{x_1}\in N(\cramped{V_1})$, and $\cramped{G_2}=G[\cramped{V_2}\cup \{\cramped{x_2}\}]$, where $\cramped{x_2}\in N(\cramped{V_2})$, and then recursively decomposing $\cramped{G_1}$ and $\cramped{G_2}$. This decomposition naturally defines an unrooted tree structure, called a \emph{split-decomposition tree}, in which the internal vertices are labeled by degenerate or prime graphs and the leaves are in bijection with $V(G)$.

The structure of the split decomposition tree of a graph may depend on the specific sequence of split operations performed on it. However, Cunningham~\cite{Cunningham82} provided a set of criteria that are satisfied by exactly one split decomposition tree for every graph. This uniqueness result is reformulated below in terms
of graph-labeled trees by Gioan~and~Paul~\cite{GiPa12}.
\begin{theorem*}[Cunningham~\cite{Cunningham82}]%
  \label{thm:cunningham}%
  For every connected graph $G$, there exists a unique split-decomposition tree such that:
  \begin{enumerate}[noitemsep, nosep]
  \item every non-leaf node has degree at least three;
  \item no tree edge links two vertices with clique labels;
  \item no tree edge links the center of a star to the extremity of
    another star.
  \end{enumerate}
\end{theorem*}
\noindent Such a tree is called \emph{reduced}, and this theorem establishes a bijection between graphs and their reduced split decomposition trees.

\section{Split-Decomposition of cactus graphs}\label{sec:characterization}

In this section, we first introduce a bijective split decomposition tree characterization of general cactus graphs. We then make certain simplifications to that characterization, which, while respecting its bijective nature, make its correspondence with our cactus grammars more apparent. Additionally, these characterizations make apparent how arbitrary size constraints on cycles of a cactus graph translate into modifications in its split-decomposition tree characterization, which can then be directly translated into modifications in the grammar.

First, we provide a characterization for the reduced split decomposition trees of general\footnote{By \emph{general}, we mean all cactus graphs with no restriction on the size of each cycle. In other words, \emph{general} cactus graphs is the family of all connected graphs in which no two cycles share an edge.} cactus graphs. The following lemma will be used in the proof of Theorem~\ref{thm:split-characterization-reduced-cactus}.

\begin{lemma}[Polygon primality]\label{lem:polygon-primality}
Polygons of size at least 5 are prime with respect to the split decomposition.
\end{lemma}
\begin{proof}
We will show that a cycle with at least 5 vertices has no splits. Suppose, on the contrary, that such a split exists.

We claim that there are at least two disjoint edges crossing the split. Any bipartition of a cycle has at least two edges crossing it. Suppose all edges crossing the bipartition are non-disjoint. Then they must all be incident to the same vertex $v$ in one side of the split. Since the degree of every vertex in a cycle is 2, there must be exactly two edges crossing the split, both of which are incident to $v$. This implies that $v$ is the only vertex in its side of the split: Any other vertex on the same side cannot have a path to $v$ without either crossing the split or causing $v$ to have degree higher than 2. Since a cycle is connected to begin with, no other vertex besides $v$ can belong to the same side of the split. This contradicts the requirement that each side of a split must have at least 2 vertices.

Let $\{x,u\}$ and $\{y,v\}$ be two disjoint edges crossing a split in a cycle with at least 5 vertices (with $x$ and $y$ on one side of the split and $u$ and $v$ on the other side). Since the edges crossing a split must induce a complete bipartite graph, the edges $\{x,v\}$ and $\{y,u\}$ must also be present in the graph. Therefore, the graph must have a cycle of size 4 as a subgraph, which is not the case for any cycle with at least 5 vertices.
\end{proof}

\begin{theorem}[\emph{Reduced} split tree characterization of general cactus graphs]%
  \label{thm:split-characterization-reduced-cactus}%
  A graph $G$ with the reduced split decomposition tree $(T,\mathcal{F})$ is a cactus graph if and only if
  \begin{enumerate}[noitemsep, nosep]
  \item $\mathcal{F}$ consists of stars and polygons of size 3 or at least 5;
  \item the center of every star-node in $T$ is attached to either
    \begin{enumerate}[label=(\alph*), noitemsep, nosep]
    \item a leaf, or
    \item the center of another star-node, only as long as both star-nodes have exactly two extremities;
    \end{enumerate}
  \item every extremity of star-nodes in $T$ is attached to either
    \begin{enumerate}[label=(\alph*), noitemsep, nosep]
    \item a polygon,
    \item a leaf, or
    \item an extremity of another star-node;
    \end{enumerate}
  \item no two polygons are adjacent.
  \end{enumerate}
\end{theorem}
\begin{proof}
  \noindent {[$\Rightarrow$]}~~Let $G$ be a cactus graph with the corresponding reduced split decomposition tree $T$. We will show that $T$ satisfies all the conditions of the theorem.
  
  A \emph{block} in a graph is a maximal 2-connected subgraph, or a subgraph formed by a bridge or an isolated vertex. In the case of a cactus graph, a \emph{block} is either a cycle, or an edge that does not belong to any cycles. Let a \emph{cluster} in a cactus graph be a maximal set of at least two blocks, all of which share exactly one vertex.
  
  Furthermore, we define \emph{external} blocks and clusters in cactus graphs in a similar way to external nodes in trees. A block in a cactus graph is called \emph{external} if it shares exactly one vertex with other blocks. A cluster is \emph{external} if it contains an external block. Note that this analogy captures the underlying tree structure for every cactus graph, where every cluster is a node in the tree, and two clusters share an edge iff they share a block. It is easy to see that the resulting graph is indeed connected and acyclic, with external clusters of the cactus graph corresponding to external nodes in the underlying tree.\footnote{However, this transformation is not bijective: There can be many cactus graphs with the same underlying tree structure. Therefore, it does not suffice for enumerative purposes.} These definitions are illustrated in Figure~\ref{fig:proof-terminology}.
  
  \begin{figure*}[ht!]
    \centering
    \includegraphics[scale=0.55]{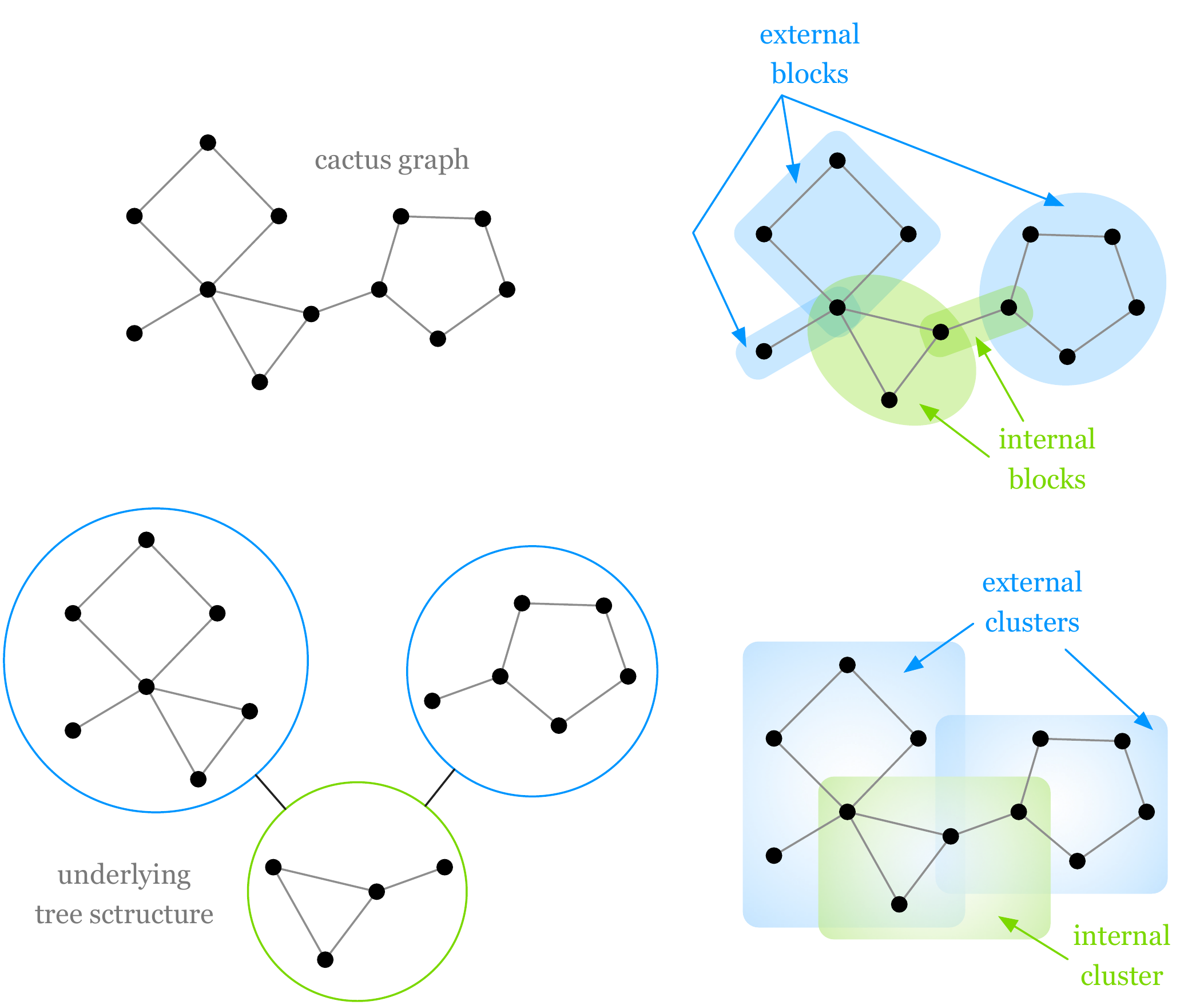}
    \caption{\label{fig:proof-terminology}%
      This figure illustrates the terminology introduced as part of the proof of Theorem~\ref{thm:split-characterization-reduced-cactus}. The top left graph is a cactus graph. The top right figure highlights all the blocks of this cactus graph, and the bottom right figure shows all its clusters. The bottom left figure shows the underlying tree structure of this cactus graph, where each node corresponds to a cluster and each edge represents a pair of clusters sharing a block. In all these figures, blue denotes \emph{external} units (\textit{i.e.} blocks, clusters, underlying tree nodes) whereas green indicates \emph{internal} units.}
  \end{figure*}

  We will proceed by induction on the number of clusters in $G$. For the base case, we consider cactus graphs with zero clusters. In this case, the graph consists of exactly one block, which can be one of the following (also shown in Figure~\ref{fig:proof-basecase}):
  \begin{itemize}[noitemsep, nosep, leftmargin=*]
      \item an edge: In this case, the corresponding reduced split decomposition tree is just an edge.
      \item a cycle of size 3 or at least 5: In this case, the corresponding reduced split decomposition tree is an internal node labeled with a cycle of the same size as the block, with a leaf hanging from every vertex of the label.~\footnote{Note that this tree is indeed reduced: A cycle of size 3 is also a clique and therefore degenerate with respect to the split decomposition, and a cycle of size at least 5 is prime with respect to the split decomposition by Lemma~\ref{lem:polygon-primality}.}
      \item a cycle of size 4: In this case, the corresponding reduced split decomposition tree is two internal nodes, each labeled with a star with two extremities. The centers of the stars are connected by an edge, and the extremities of the stars are connected to leaves.
  \end{itemize}
  It is easy to confirm that the split decomposition trees above are both correct (by making sure their alternated paths match the adjacencies of the corresponding cluster) and reduced (by making sure the conditions of Theorem~\ref{thm:cunningham} are met). Furthermore, all the above split decomposition trees satisfy the conditions of this theorem, thus proving our base case.
  
  \begin{figure*}[ht!]
    \centering
    \includegraphics[scale=0.55]{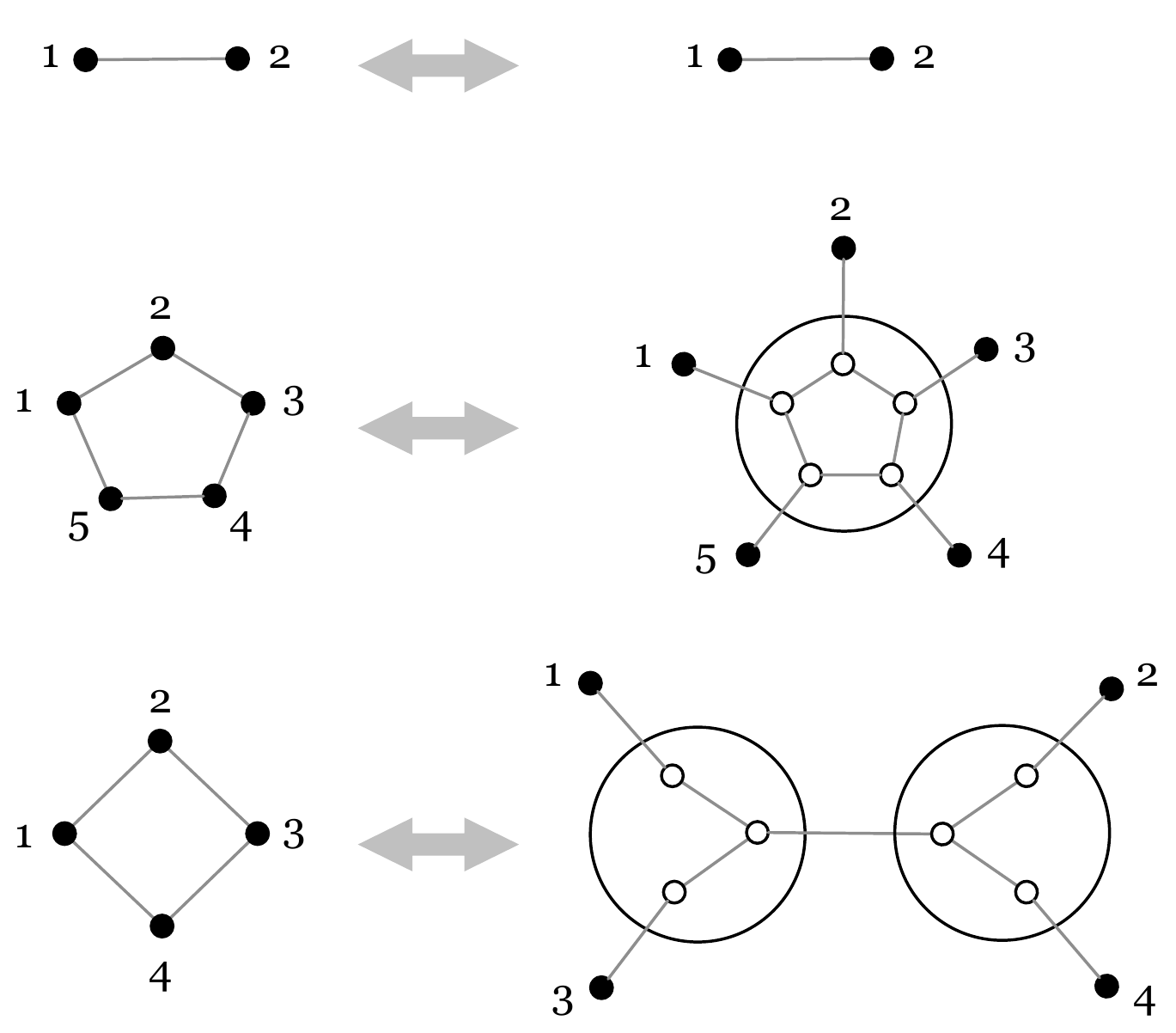}
    \caption{\label{fig:proof-basecase}%
      This figure shows the base case in part of the proof of Theorem~\ref{thm:split-characterization-reduced-cactus}, \textit{i.e.} split decomposition trees of cactus graphs with exactly one block (and thus zero clusters). These graphs are partitioned into three disjoint subsets: an edge, a cycle of size 3 or at least 5, and a cycle of size 4. Each subset is treated differently as illustrated above. Note that the numbers here serve to specify the bijection between the vertices of each cactus graph and the leaves of its split decomposition tree. They do not serve to distinguish the vertices and do not imply that these graphs are labeled.}
  \end{figure*}
  \noindent
      
  For the inductive step, let $G$ be a cactus graph with at least one cluster. First, we claim that $G$ must have an external cluster, since the underlying tree must have a leaf. Let $C$ be an external cluster including blocks $\cramped{S_1}\cdots\cramped{S_k}, k>1$, all of which share the vertex $v\in V(G)$. Finally, observe that at most one of the blocks is internal (the block representing the edge between the node corresponding to $C$ and its neighbor in the underlying tree is internal). Without loss of generality, assume that $\cramped{S_k}$ is the internal block, if one exists.
  
  We split the cluster $C$ from the rest of the cactus graph and create a graph-labeled tree $\tilde{T}$ as follows (Figure~\ref{fig:proof-step}):
  \begin{itemize}[noitemsep, nosep, leftmargin=*]
      \item Create a star-node $s$ with $k$ extremities, one per block, respecting the blocks' circular permutation around $v$ in the case of plane cacti. Attach the center of this star-node to a leaf corresponding to $v$.
      \item For every external block $\cramped{S_i}, i<k$,
      \begin{itemize}[noitemsep, nosep, leftmargin=*]
        \item If $\cramped{S_i}$ is a single edge, attach its corresponding extremity to a leaf.
        \item If $S_i$ is a cycle of size 4, create two star-nodes with two extremities each, where the pair of extremities of each star represent two non-adjacent vertices in $\cramped{S_i}$, with the order of the appearance of the extremities respecting the planar embedding of the corresponding vertices in the case of plane cacti. Attach the two stars together at their centers, and connect the extremity corresponding to $v$ with the extremity of $s$ corresponding to $\cramped{S_i}$, and connect the other three extremities to leaves.
        \item Otherwise, create an internal node labeled with $\cramped{S_i}$. Attach the copy of $v$ in the label of the new internal node to the extremity of $s$ corresponding to $\cramped{S_i}$. Attach every other vertex in the label of the new internal node to a leaf, representing the vertices of $G$ in $\cramped{S_i}$, respecting their planar ordering.
      \end{itemize} 
      \item Next, consider $G'$, the remainder of $G$ excluding the external blocks $\cramped{S_i},i<k$ but including $\cramped{S_k}$. This is a cactus graph with fewer clusters than $G$. Applying the inductive hypothesis, we obtain a reduced split-decomposition tree $\tilde{T}'$ satisfying the conditions of this theorem. Finally, we attach $\tilde{T}'$ to $s$ by removing the leaf corresponding to $v$ in $V(\tilde{T}')$ and connecting the edge incident to that leaf with the extremity of $s$ corresponding to $\cramped{S_k}$.
  \end{itemize}
  Note that this procedure produces a split-decomposition tree $\tilde{T}$ that satisfies the conditions of this theorem. In every step, the new star-node and the nodes corresponding to external blocks satisfy the conditions by construction, and the nodes within $G'$ satisfy the conditions by the inductive hypothesis. It remains to check that nothing is violated upon attaching $\tilde{T}'$ to an extremity of the new star-node. Since $v$ belongs to a single external block in $G'$, its corresponding leaf in $\tilde{T}'$ can only be attached to either a star extremity or a polygon; this is because the only time a leaf is attached to a star center is when the leaf corresponded to the common vertex of a cluster in the inductive step, in which case it was shared between more than one block. Therefore, only conditions 3 and 4 are relevant, both of which are satisfied when attaching the extremity of $s$ to $v$'s neighbor in $\tilde{T}'$.
  
  \begin{figure*}[ht!]
    \centering
    \includegraphics[scale=0.55]{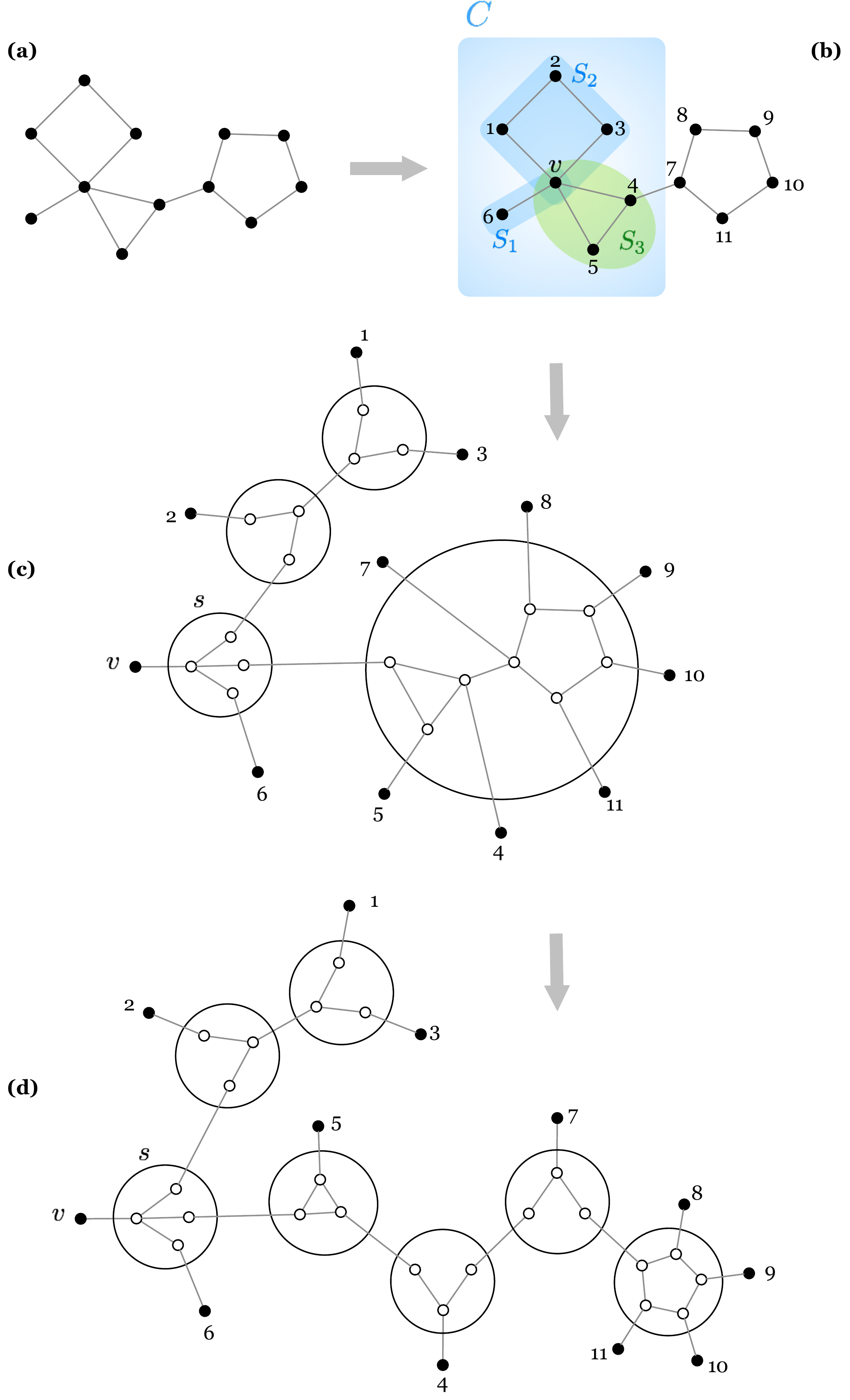}
    \caption{\label{fig:proof-step}%
    This figure shows the inductive step used in the proof of Theorem~\ref{thm:split-characterization-reduced-cactus}, \textit{i.e.}. The inductive step begins by identifying an external cluster $C$, which consists of blocks $\cramped{S_1}\cdots\cramped{S_k}$, all of which share a vertex $v$. The external blocks (here $\cramped{S_1}$ and $\cramped{S_2}$) are then separated via a star node $s$, the center of which is attached to $v$. The rest of the cactus graph forms a large super node hanging from $s$, and will be decomposed using the inductive hypothesis. Furthermore, the specific plane embedding of the original cactus graph around $v$ can be respected by the choice of the order of the extremities of $s$. Similar to Figure~\ref{fig:proof-basecase}, the numbers here are merely for illustrative purposes and do not imply that the cactus graph is labeled.}
  \end{figure*}
  \noindent
  
  Finally, we need to show that this procedure preserves the adjacencies of $G$.
  
  First, we show that every step in the procedure above preserves the adjacencies of $G$ via alternated paths in $\tilde{T}$. At every step, any edge $e=\{x,y\}\in E(G)$ must belong to one of the following cases.
  \begin{itemize}[noitemsep, nosep, leftmargin=*]
      \item $e$ is incident to $v$: In this case, $e$ must belong to a block $\cramped{S_i}, i\leqslant k$. There is an alternated path from $v$ to the center of the star-node attached to $v$, then to the extremity corresponding to $\cramped{S_i}$ and then out of the star-node. From there,
      \begin{itemize}[noitemsep, nosep, leftmargin=*]
          \item if $\cramped{S_i}$ is an edge, we arrive at a leaf representing the other end of $e$, completing our alternated path;
          \item if $\cramped{S_i}$ is an external cycle of size 4, we next arrive at the extremity of a star-node and can continue to its center, then to the center of the neighboring star-node, then to the extremity of the current star-node corresponding to the other end of $e$, and finally to the leaf corresponding to the other end of $e$;
          \item finally, if $\cramped{S_i}$ is an external cycle of size 3 or at least 5 or an internal cycle, we next arrive at an internal node labeled with either $\cramped{S_i}$ or $G'$, both of which include $e$ as an edge within the label. We can therefore continue the alternated path by taking the edge corresponding to $e$ from within the label and then to the leaf corresponding to the other end of $e$.
      \end{itemize}
      \item $e$ is not incident to $v$: In this case, $e$ belongs to either $\cramped{S_i}, i<k$ or $G'$. Note that in the first case $\cramped{S_i}$ cannot be an edge since $v\not\in e$.
      \begin{itemize}[noitemsep, nosep, leftmargin=*]
          \item If $\cramped{S_i}$ is a cycle of size 4, there is an alternated path representing $e$ starting from the leaf corresponding to $x$ to an extremity of the neighboring star, then to its center, then to the center of the next star, and finally to the extremity of the same star corresponding to the leaf $y$ and lastly to the leaf corresponding to $y$ itself.
          \item Otherwise, $\cramped{S_i}$ is either cycle of size 3 or at least 5, or $G'$. In all these cases, $e$ is present as within the label of an internal node, and thus there is an alternated path between the leaves representing $x$ and $y$ that takes that edge within a graph label corresponding to $e$.
      \end{itemize}
  \end{itemize}
  Therefore, at every step, including the base case, for every edge in $G$ there exists an alternated path in $\tilde{T}$.
  
  In the other direction, we have to show that in every step, for every alternated path in $\tilde{T}$ there is an edge in $G$. Take a pair of vertices $x,y\in V(G), \{x,y\}\not\in E(G)$. It is easy to check that if $x$ and $y$ belong to the same block $\cramped{S_i}$ of size at least 5, any path between the leaves representing them in $\tilde{T}$ uses more than one edge from the the internal node labeled with $\cramped{S_i}$. Furthermore, if $x$ and $y$ belong to the same block of size 4 or different blocks, any path connected the leaves representing them uses two edges from the same star-node.
  
  Finally, note that the conditions of this theorem imply that $T$ is \emph{reduced}: Cunningham's theorem's (Section~\ref{thm:cunningham}) criterion 1 for being reduced is satisfied by condition 1 of this theorem, criterion 2 by condition 4 of this theorem, and criterion 3 by conditions 2 and 3 of this theorem. Therefore, the procedure above produces a reduced split-decomposition tree, which is guaranteed to be unique by Cunningham's theorem. Thus, $\tilde{T}$ and $T$ are the same split-decomposition tree.
  
  Therefore, the adjacencies of $G$ correspond bijectively to the alternated paths of $T$. Furthermore, the plane embedding of a cactus graph can be preserved if the cyclic order of the new star-node in every step respects the plane embedding of the corresponding blocks in $G$.
  
  % % % % % % % % % % % % % % % % % % % % % % % % % % % % % % %
  
  \noindent {[$\Leftarrow$]}~~Let $T$ be a reduced split-decomposition tree satisfying the conditions of this theorem. We should show that the original graph corresponding to $T$ is indeed a cactus.
  
  In the first part of the proof, we introduced a procedure for constructing the reduced split-decomposition tree of any cactus graph. The reverse of the same procedure can be used to build the original cactus graph from a split-decomposition tree. More specifically, pairs of star-nodes connected via their centers can be combined to form cycles of size 4, and star-nodes with a leaf hanging from their center can be combined with their neighbors in the reverse direction of the procedure as well. It is easy to check that as long as $T$ has more than one internal node, it is possible to combine some internal nodes into one. Furthermore, these reverse steps maintain the invariant that in all intermediate split-decomposition trees, all the created labels are cactus graphs. Therefore, at the end of this procedure, a split decomposition tree with a single internal node remains, the label of which must also be a cactus graph.
  % % % % % % % % % % % % % % % % % % % % % % % % % % % % % % %
\end{proof}

Theorem~\ref{thm:split-characterization-reduced-cactus} gives a bijective characterization of cactus graphs in terms of \emph{reduced} split-decomposition trees. In its proof we utilize Cunningham's Theorem, which guarantees uniqueness of the reduced split-decomposition tree for any graph. It turns out, however, that the characterization for general cactus graphs can be simplified to the form in Theorem~\ref{thm:split-characterization-cactus}, using a bijection between the reduced split-decomposition trees in Theorem~\ref{thm:split-characterization-reduced-cactus} and a simpler, yet not reduced, set of split-decomposition trees.

\begin{theorem}[Split tree characterization of general cactus graphs]%
  \label{thm:split-characterization-cactus}%
  Cactus graphs are in bijection with split-decomposition trees where
  \begin{enumerate}[label=(\alph*), noitemsep, nosep]
  \item internal nodes are stars and polygons;
  \item the centers of all star-nodes are attached to leaves;
  \item the extremities of star-nodes are attached to leaves or polygons.
  \end{enumerate}
\end{theorem}
\begin{proof}
Since the characterization in Theorem~\ref{thm:split-characterization-reduced-cactus} describes a \emph{reduced} tree, a direct application of Cunningham's Theorem guarantees a one-to-one correspondence between cactus graphs and those split decomposition trees. It suffices to show a bijection between the characterization given in Theorem~\ref{thm:split-characterization-reduced-cactus} and the characterization in this Theorem. This is easily achieved by combining every pair of star-nodes with their centers attached into cycles of size 4, while leaving everything else intact.
\end{proof}

%%%%%%%%%%%%%%%%%%%%%%%%%%%%%%%%%%%%%%%%%%%%%%%%%%%%%%%%%%%%%%%%%%%%%%%%%%

\begin{table*}
  \centering
  \def\arraystretch{1.2}
  \begin{tabularx}{0.84\textwidth}{cW}
    \toprule
    \textbf{Symbol}&\textbf{Explanation}\\
    \midrule
    %%%%%%%%%%%%%%%%%%%%%%%%%%%%%%%%%%%%%%%%%%%%%%%%%%%%%%%%%%%%%%%%%%%%%%
    $\cls{P}$ &{a polygon-node entered from one of its vertices (and missing the corresponding subtree)}\\
    $\cls[C]{S}$ &{a star-node entered through its \emph{center} (and missing the corresponding subtree)}\\
    $\cls[X]{S}$ &{a star-node entered through one of its (at least two) \emph{extremities} (and missing the corresponding subtree)}\\[0.4em]
    %%%%%%%%%%%%%%%%%%%%%%%%%%%%%%%%%%%%%%%%%%%%%%%%%%%%%%%%%%%%%%%%%%%%%%
    $\clsAtom$ &{a leaf of the split-decomposition tree
                 (an atom with unit size)}\\[0.4em]
    %%%%%%%%%%%%%%%%%%%%%%%%%%%%%%%%%%%%%%%%%%%%%%%%%%%%%%%%%%%%%%%%%%%%%%
    \midrule
    %%%%%%%%%%%%%%%%%%%%%%%%%%%%%%%%%%%%%%%%%%%%%%%%%%%%%%%%%%%%%%%%%%%%%% 
    $\cls[P]{T}$ &{a {split-decomposition tree} rerooted at a
                   polygon-node (all subtrees are present)}\\
    $\cls[S]{T}$ &{a {split-decomposition tree} rerooted at a
                   star-node (all subtrees are present)}\\[0.4em]
    %%%%%%%%%%%%%%%%%%%%%%%%%%%%%%%%%%%%%%%%%%%%%%%%%%%%%%%%%%%%%%%%%%%%%% 
    $\cls[SP]{T}$ &{a {split-decomposition tree} rerooted at an
                     \emph{edge} connecting a star-node to a polygon-node (the edge can either connect the polygon-node to the star-node's center or an extremity; the edge accounts for one subtree of the polygon-node and one subtree of the star-node)}\\[0.4em]
    %%%%%%%%%%%%%%%%%%%%%%%%%%%%%%%%%%%%%%%%%%%%%%%%%%%%%%%%%%%%%%%%%%%%%% 
    \bottomrule
  \end{tabularx}
  \caption{\label{tab:symbols}%
    The main symbols used to define the {split-decomposition tree} of the cactus graphs. Refer to Subsection~\ref{subs:split-grammars} for details on the terminology, and to prior papers~\cite[\S 1.6]{BaLu17} for details on the dissymmetry theorem, from which all the rerooted trees, denoted by $\cramped{\cls[\omega]{T}}$, come from.}
\end{table*}

%%%%%%%%%%%%%%%%%%%%%%%%%%%%%%%%%%%%%%%%%%%%%%%%%%%%%%%%%%%%%%%%%%%%%%%%%%

\section{Generating cactus grammars}\label{sec:grammars}
In this section we provide guidelines on how to produce grammars\footnote{Grammars are a standard concept in computer science, since they are used to specify programming languages. Although we assume the semantics of these grammars is straightforward, we briefly recall the various operators in Subsection~\ref{subs:FS}.} for different varieties of cactus graphs. These guidelines do not assume familiarity with the technical details of the proofs in the previous sections and can be applied as a black box. The following subsections describe how a grammar can be obtained from a a characterization of a split-decomposition tree, such as the one given in Theorem~\ref{thm:split-characterization-cactus}. For the black box version of the grammars, the reader can skip to Section~\ref{subs:template}.

\subsection{Decomposable structures.\label{subs:FS}}

In order to enumerate classes of split-decomposition trees, we use the framework of decomposable structures, described by Flajolet and Sedgewick~\cite{FlSe09}. We refer the reader to this book for details and outline below the basic idea.

We denote by $\clsAtom$ the combinatorial family composed of a single object of size $1$, usually called \emph{atom} (in our case, these refer to a leaf of a split-decomposition tree, \textit{i.e.}, a vertex of the corresponding graph).

Given two disjoint families $\cls{A}$ and $\cls{B}$ of combinatorial objects, we denote by $\cls{A} + \cls{B}$ the \emph{disjoint union} of the two families and by $\cls{A} \times \cls{B}$ the \emph{Cartesian product} of the two families.

We denote by $\Set{\cls{A}}$ (resp. $\cramped{\Set[\Omega]{\cls{A}}}$ where $\Omega\subseteq{\cramped{\mathbb{Z}_{\geqslant0}}}$) the family defined as all sets (resp. sets of sizes that belong to $\Omega$) of objects from ${\cls{A}}$\footnote{Note the $\Set$ operator refers to multisets in the unlabeled world and sets in the labeled world.}.  In a similar fashion, we define $\Seq{\cls{A}}$ (resp. $\cramped{\Seq[\Omega]{\cls{A}}}$) as the family of all (ordered) sequences of objects from ${\cls{A}}$ (resp. a sequence the number of objects in which belongs to $\Omega$). Finally, we define cycles $\Cyc{\cls{A}}$, undirected sequences $\Useq{\cls{A}}$, and undirected cycles $\Ucyc{\cls{A}}$ similarly.

\TODOAll{Add details about the constructions}

\subsection{Split-decomposition trees expressed symbolically.%
  \label{subs:split-grammars}}

While approaching graph enumeration from the perspective of tree
decomposition is not a new idea (the recursively decomposable nature of trees makes them well suited to enumeration), Chauve~\etal~\cite{ChFuLu17} brought specific focus to Cunningham's split-decomposition. The same framework was used by Bahrani and Lumbroso~\cite{BaLu17} and is the starting point of this paper, so we briefly outline their method here.

\paragraph{Example} Let us consider the split-decomposition tree drawn in Figure~\ref{fig:proof-step} (d), and illustrate how it can be expressed recursively as a rooted tree.

Suppose the tree is rooted at vertex $v$. Assigning a root immediately defines a direction for all tree edges, which can be thought of as oriented away from the root. Starting from the root, we can set out to traverse the tree in the direction of the edges, one internal node at a time.

We start at the root, vertex $v$. The first internal node we encounter is a star node, $s$, and since we are entering it from the star's center, we have to describe what is on each of its three extremities.

\begin{remark*}
If we are considering a \emph{plane} cactus graph, the cyclic permutation of the edges incident to $v$ in the original graph defines a cyclic permutation on the extremities of $s$, so these extremities can be described as a \emph{cycle}. On the other hand, if we are considering a \emph{free} cactus graph that is not embedded in the plane, there is no imposed ordering on the extremities of $s$. In that sense, the extremities of $s$ are indistinguishable and can therefore be described as a \emph{set}. These differences are what differentiates the templates in Section~\ref{subs:template}, more specifically the first row of Table~\ref{tab:operator-guidelines}.
\end{remark*}

Next, we will visit the extremities of $s$ (in counter-clockwise order for easier referencing). On the first extremity, there is a leaf, 6, which concludes our journey in that direction; on each of the other two extremities, there is another split-decomposition subtree, which we can continue to explore.

The second extremity of $s$ leads to a split-decomposition subtree, of which the first internal node we encounter is labeled by a polygon of size 3. Continuing the traversal through this polygon, we need to describe what hangs from each vertex of the polygon. The exploration continues recursively through all the split-decomposition subtrees hanging from the vertices of this polygon until leaves are reached.
\begin{remark*}
Once again, there is a distinction between the plane and non-plane cases. In the case of a \emph{plane} cactus graph, each vertex of this polygon comes with a cyclic permutation of its edges, which uniquely determines the ``next'' vertex in the polygon, thus imposing a direction with respect to which the polygon can be traversed. Therefore, the split-decomposition subtrees hanging from the remaining vertices of the current polygon can be described as a \emph{sequence}. On the other hand, in the case of a \emph{free} cactus graph, there is no such imposed direction, and the split-decomposition subtrees can be described as an \emph{undirected sequence}.
\end{remark*}
\begin{remark*}
Note that unlike the extremities of a star, which are completely indistinguishable in the non-plane setting, the vertices of a polygon come with a cyclic ordering defined by the edges regardless of plane embedding, and therefore form an undirected sequence rather than a set. This distinction is reflected in the third row of Table~\ref{tab:operator-guidelines}.
\end{remark*}

The split-decomposition subtree hanging from the third extremity of $s$ is traversed recursively in a similar manner, covering the whole tree.

\paragraph{Grammar description} The characterization developed in Theorem~\ref{thm:split-characterization-cactus} can be translated to a symbolic grammar following the same approach as outlined in the example above. As an example, consider the class of plane, rooted, pure 5-cacti. To describe the split-decomposition tree corresponding to this class symbolically, we must give a set of rules that prescribe the possible neighbors for each type of internal node.

Let's consider the rule for star nodes. First, assume like at the beginning of our example, that we enter a star node through its center; we have to describe what the extremities can be connected to.

According to Cunningham's Theorem, we know that there are at least two extremities (since every non-leaf node has degree at least three). Furthermore, by Theorem~\ref{thm:split-characterization-cactus}, we know that star extremities can be connected to either leaves or polygons. More specifically, connecting a star extremity to a leaf corresponds to a pendant edge in the original graph, whereas connecting a star extremity to an internal node labeled by a polygon of size $k$ corresponds to a polygon of that size in the original graph. We call $\cramped{\cls[C]{S}}$ a split-decomposition tree that is traversed starting at a star-node entered through its center. We have
\begin{align*}
  \cls[C]{S} = \Cyc[\geqslant 2]{\cls{P}}
\end{align*}
because indeed, we have at least two extremities, which have a cyclic order since we are describing a \emph{plane} cactus graph, and each of these extremities must lead to an internal node labeled by a polygon, $\cls{P}$.

A split decomposition subtree rooted at an internal node labeled by a polygon $\cls{P}$ entered from one of its vertices can in turn be described by a symbolic equation
\begin{align*}
  \cls{P} = \Seq[4](\clsAtom + \cls[X]{S}),
\end{align*}
where $\clsAtom$ denotes a leaf (an atom of size 1) and  $\cls[X]{S}$ denotes a star node entered through one of its extremities. Indeed, since we are in a \emph{plane} setting, the split-decomposition subtrees hanging from the remaining vertices of the current polygon can be described as a \emph{sequence}]. The cardinality constraint captures the \emph{pure} aspect of these cacti. Finally, by Theorem~\ref{thm:split-characterization-cactus}, internal nodes with polygonal labels can be attached to leaves or star extremities, which translates into the disjoint union term $\clsAtom + \cls[X]{S}$.

The rule for a star node entered through its extremity is derived in a similar manner.

\paragraph{Conventions} As explained above, we use rather similar
notations to describe the combinatorial classes that arise from
decomposing split-decomposition trees. These notations are summarized in Table~\ref{tab:symbols}.

\paragraph{Terminology} In the rest of this paper, we describe the combinatorial class $\cramped{\cls[X]{S}}$ as representing a ``\emph{a star-node entered through an extremity}'', but others may have alternate descriptions: such as ``\emph{a star-node linked to its parent by an extremity}''; or such as Iriza~\cite{Iriza15}, ``\emph{a star-node with the subtree incident to one of its extremities having been removed}''---all these descriptions are equivalent (but follow different viewpoints).

%%%%%%%%%%%%%%%%%%%%%%%%%%%%%%%%%%%%%%%%%%%%%%%%%%%%%%%%%%%%%%%%%%%%%%%%%%

\subsection{General template.}\label{subs:template}

The variations of cactus graphs supported by our grammars include \textit{plane} versus \textit{non-plane} and \textit{rooted} versus \textit{unrooted}. These differences gives rise to different operators in the grammars, which we will outline below.

Another variation in types of cacti is whether the vertices are labeled. It turns out \textit{labeled} versus \textit{unlabeled} cacti have identical grammars but are distinguished by the way their grammars are translated to generating function equations, with exponential generating functions for labeled cacti and ordinary generating functions for unlabeled ones.

A generalization of \textit{pure} versus \textit{mixed} cacti is supported in the form of arbitrary constraints on allowable cycle sizes. More specifically, in our grammars $\Omega$ refers to a set of integers each greater than 1, such that the size of each cycle in the cactus graph must belong to that set. For example, in the case of \textit{pure} cacti, $\Omega$ contains a single element, and setting $\Omega=\set{2}$ gives the grammar for trees. Furthermore, the grammar for \textit{mixed} cacti can be obtained by letting $\Omega$ be the set of all integers greater than 1. Finally, we use the notation $\Omega-1\coloneqq\{x-1\text{ for all } x\in\Omega\}$.

The grammars for all the various cactus graphs mentioned above roughly follow the same template. We outline this template below, highlighting the minor differences from case to case.
The following three equations are common to all grammars, with $\Xxx$, $\Yyy$, and $\Zzz$ replaced by the appropriate operator as listed in Table~\ref{tab:operator-guidelines}.
\begin{align*}
    \cls[C]{S}  &= \Xxx[\geqslant2]{\cls{P}}\\
    \cls[X]{S}  &= \clsAtom\times\Yyy[\geqslant1]{\cls{P}}\\
    \cls{P}     &= \Zzz[\Omega-1]{\clsAtom+\cls[X]{S}}
\end{align*}
The class of rooted cactus graphs $\cls[\mLeaf]{G}$ can be specified as
\begin{align*}
    \cls[\mLeaf]{G}  = \clsAtom\times(\cls{P}+\cls[C]{S}).
\end{align*}
Furthermore, the class of unrooted cactus graphs, $\cls{G}$, can be specified as
\begin{align*}
    \cls{G}     &= \cls[S]{T} + \cls[P]{T} - \cls[SP]{T}\\
    \cls[S]{T}  &= \clsAtom\times\cls[C]{S}\\
    \cls[P]{T}  &= \Www[\Omega]{\clsAtom+\cls[X]{S}} \\
    \cls[SP]{T} &= \cls{P}\times\cls[X]{S},
\end{align*}
with $\Www$ replaced by the appropriate operators according to Table~\ref{tab:operator-guidelines}.

\begin{table}[h!]
\centering
\begin{tabular}{ c|c|c } 
       & \textbf{plane} & \textbf{free} \\
 \hline
 $\Xxx$ & \Cyc & \Set \\ 
 $\Yyy$ & \Seq & \Set \\ 
 $\Zzz$ & \Seq & \Useq \\
 $\Www$ & \Cyc & \Ucyc
\end{tabular}
\caption{Guidelines for choosing operators in the template grammar for \textit{plane} versus \textit{free} cactus graphs.}
\label{tab:operator-guidelines}
\end{table}

\paragraph{Unrooting the grammar from the traversal of the tree}
The grammar template $\cls[\mLeaf]{G}$ corresponds to \emph{rooted} cactus graphs, and results directly from applying the constraints on the split-decomposition tree---which we defined in Section~\ref{sec:characterization}---by using the same type of process as described earlier in Subsection~\ref{subs:split-grammars}. But obtaining grammars for the \emph{unrooted} cactus graphs requires additional tools, as there are symmetries to account for, that depend on the structure of individual graphs.
    
There are at least two such \emph{unrooting} techniques. We apply here the \emph{dissymmetry theorem [for trees]}, from Bergeron \textit{et~al.}~\cite{BeLaLe98}, with a proof in Drmota~\cite[\S 4.3.3, p.~293]{HoECDrmota15}. In the context of symbolically specified split-decomposition trees, Chauve~\etal~\cite[\S 2.2 and \S 3]{ChFuLu17} contain indications on how to apply the theorem concretely, with further remarks by Bahrani~\etal~\cite[\S 1.6]{BaLu17}.
    
The dissymmetry theorem is relatively easy to apply, but the equation which has a subtraction is algebraic but not symbolic: the subtraction is meaningful when we are tabulating the enumeration, but it has no combinatorial interpretation. This means such grammars cannot benefit from the direct application of some of the most powerful ``last kilometer'' results from analytic combinatorics transfer theorems: asymptotics, parameter analysis, random generation (and so empirical analysis).
    
An alternate tool to unroot combinatorial classes,  \emph{cycle-pointing}~\cite{BoFuKaVi11}, does not have this issue. It is a combinatorial operation (rather than an algebraic one), and preserves the transfer theorems, in particular the results on random generation---the recursive method~\cite{NiWi75, FlZiVa94, DeZi99, PoTeDe06} and Boltzmann sampling \cite{DuFlLoSc04, PiSaSo12, BoLuRo15}. It involves a more complex decomposition of the grammar and finer understanding of the symmetries of the structure. Iriza~\cite{Iriza15} has applied it to the distance-hereditary and 3-leaf power grammars of Chauve~\etal~\cite{ChFuLu17}. Iriza's work was the basis of the random samplers used for Figure~\ref{fig:random-cacti}.

\subsection{Example.}

\TODOAll{Simpler Example}

\begin{figure*}[ht!]
    \centering
    \includegraphics[scale=0.6]{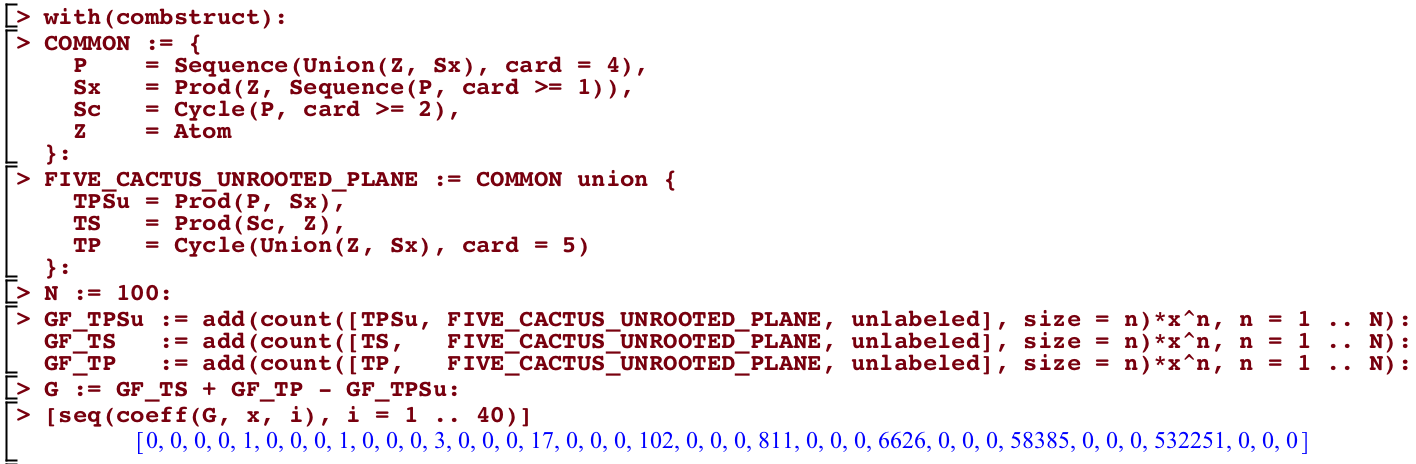}
    \caption{\label{fig:5cactus-maple}%
      The enumeration of different varieties of cactus graphs can be derived using our grammars and Maple. This figure shows a Maple worksheet for computing the enumeration of unrooted, unlabeled, pure, plane 5-cacti.}
\end{figure*}

We state a grammar for plane unrooted pure 5-cacti ($\cls[\mathrm{pu5}]{G}$) as an example:
\begin{align*}
    \cls[\mathrm{pu5}]{G}  &= \cls[S]{T} + \cls[P]{T} - \cls[SP]{T}\\
    \cls[S]{T}  &= \clsAtom\times\cls[C]{S}\\
    \cls[P]{T}  &= \Cyc[\{5\}]{\clsAtom+\cls[X]{S}} \\
    \cls[SP]{T} &= \cls{P}\times\cls[X]{S}\\
    \cls[C]{S}  &= \Cyc[\geqslant2]{\cls{P}}\\
    \cls[X]{S}  &= \clsAtom\times\Seq[\geqslant1]{\cls{P}}\\
    \cls{P}     &= \Seq[\{4\}]{\clsAtom+\cls[X]{S}}.
\end{align*}
\noindent
A symbolic math package can be used to extract the exact enumeration of plane unrooted pure 5-cacti from this grammar in polynomial time (See Figure~\ref{fig:5cactus-maple} for a sample implementation in Maple). The resulting counting sequence in this case indeed matches the sequence provided by Bona~\textit{et al.}~\cite{bona2000enumeration}. The first few terms of the enumeration, \EIS{A054365}, are listed below.
\begin{align*}
	&1, 1, 1, 3, 17, 102, 811, 6626, 58385,\\
	&\qquad 532251, 5011934, 48344880, \dots\text{.}
\end{align*}

\TODOAll{Reviewer comment: As the paper is on the enumeration of cactus graphs, I would have expected integer sequences for cactus graphs (for example as given in Section 1.2 on page 3). How do the specifications in Sec 1.2 relate to the Husimi trees as given in the OEIS? How do your results relate to approaches used in Bergeron et al. (Chapter 4.2, and specifically page 301 (Example 1)).}

%%\TODOAll{include the typeset maple file}

%%\TODOJeremie{Try tweaking Sam's code to get UCYC and USEQ working}

\subsection{Simplified grammars.} 

In this section, we list simplified versions of the grammars produced by the template in the previous section. The previous template is informative in its maintaining of a general common structure across different sub-classes of cactus graphs. More specifically, in the previous template, the differences between grammars of \textit{similar} classes are relatively small. The following simplified grammars, however, are more concise and therefore amenable to further analysis, such as the extraction of asymptotics and the study of parameters.

\begin{itemize}
    \item \textit{free, rooted} cactus graphs (denoted by $\clsFRG$) admit the following recursive grammar:
    \begin{align*}
        \clsFRG = \clsAtom \times  \Set[\geqslant1]{\Useq[\Omega-1]{\clsAtom+\clsFRG}}
    \end{align*}
    \item \textit{plane, rooted} cactus graphs (denoted by $\clsPRG$) can be specified as
    \begin{align*}
        \clsPRG &= \clsAtom \times  \Cyc[\geqslant1]{\Seq[\Omega-1]{\cls{Q}}}\\
        \cls{Q} &= \clsAtom\times\Seq{\Seq[\Omega-1]{\cls{Q}}}
    \end{align*}
    \item \textit{free, unrooted} cactus graphs (denoted by $\clsFUG$) can be specified as
    \begin{align*}
        \clsFUG &= \Ucyc[\Omega]{\cls{Q}} - \cls{Q} \times\Useq[\Omega-1]{\cls{Q}} + \cls{Q}-\clsAtom\\
        \cls{Q} &= \clsAtom\times\Set{\Useq[\Omega-1]{\cls{Q}}}
    \end{align*}
    \item \textit{plane, unrooted} cactus graphs (denoted by $\clsPUG$) can be specified as
    \begin{align*}
        \clsPUG &= \Cyc[\Omega]{\cls{Q}} + \clsAtom\times\Cyc[\geqslant2]{\Seq[\Omega-1]{Q}}\\
        &\ \ \ \ - (\cls{Q}-\clsAtom)\times\Seq[\Omega-1]{\cls{Q}}\\
        \cls{Q} &= \clsAtom\times\Seq{\Seq[\Omega-1]{\cls{Q}}}
    \end{align*}
\end{itemize}

\noindent With these grammars, and using existing frameworks and tools from analytic combinatorics, many results on cactus graphs are easily accessible: full exact enumeration (see Figure~\ref{fig:5cactus-maple}), asymptotics, and exhaustive and random generation~\cite{DuFlLoSc04}.

%%%%%%%%%%%%%%%%%%%%%%%%%%
%%% SECTION Conclusion %%%

\section{Conclusion}

This paper continues the line of work on split-decomposition from an analytic combinatorial perspective, initiated by Chauve \textit{et~al.}~\cite{ChFuLu14, ChFuLu17} and which we have also explored previously, in the context of subclasses of distance-hereditary graphs defined by forbidden subgraphs~\cite{BaLu17}. While these works focused on graphs that were fully decomposable by the split-decomposition, we tentatively investigate having a certain type of prime nodes: cycles.\footnote{Even though this class of prime nodes may seem trivial to consider, consider some other examples of graph class that decomposes with a well-defined set of prime nodes. For instance, \emph{circle graphs}~\cite{Spinrad94} have a subset of circle graphs as prime nodes, therefore likely requiring some careful substitution and mutual recursion. Another example are parity graphs, which admit as prime node any bipartite graph; we have known how to tabulate the enumeration of bipartite graph for a long time~\cite{Hanlon79}, but it is only recently that Gainer-Dewar and Gessel~\cite{GaGe14}  exactly enumerated bipartite graphs with some preservation of the symmetries, but using sophisticated algebraic tools that do not seem to translate easily to symbolic specifications, or at least not at first glance.}

In this paper, we looked at the split-decomposition of cactus graphs, and provided a characterization of the corresponding split-decomposition trees, as Gioan and Paul~\cite{GiPa12} did for cographs, 3-leaf power graphs and distance-hereditary graphs. We then leveraged this characterization, as Chauve~\etal~\cite{ChFuLu14, ChFuLu17} did, to obtain a symbolic specification for cactus graphs---and using the same too~\cite{Drmota97}. 

Furthermore, because cactus graphs are akin to trees, of which the nodes have been replaced by cycles, the structure of their split-decomposition tree is very similar to that of the cactus graph it decomposes. It was thus possible to define a convenient template of grammars to symbolically describe any class of rooted cactus graph.

This offers perhaps the possibility of more systematic results on cactus graphs. Perhaps some general theorem on the asymptotic that depends only on the allowed size of cycle. We also hope that the convenient format of the grammars will make certain other applications from analytic combinatorics more accessible. For instance, the random graphs in Figure~\ref{fig:random-cacti} were produced by our rooted cactus graph sampler.

Along the lines of the split decomposition, it would be interesting to consider more general graphs and split-decomposition trees with other classes of prime nodes, especially in the case of families of graphs that are less "tree-like" than cacti. A good starting point can be considering split-decomposition trees where the labels are stars, cliques, or prime cycles (a superset of distance-hereditary graphs).

Another approach that might hold some interest: define and enumerate graph classes according to the subset of allowable prime nodes, from among specific graph classes---though it is not clear which graph classes those could be.

\section*{Acknowledgements}
We would like to thank the anonymous reviewers for critically reading the manuscript and making suggestions that led to significant improvements to the content and clarity of the paper.
% \TODOAll{mention forbidden paper appendix on 4cacti and how our grammars are the same + combinatorial explanation.} -> Maybe for next version

%\section{What can be done with grammars}
%\TODOAll{sample enumerations. mention what enumerations were unknown. mention generation. remark on non-generativity of unrooted grammars and future work. Add remark on how having the unrooted grammar is important in the unlabeled case for enumeration since can't just divide by n}
\newpage
\TODOMaryam{Add methodology chain figure from talk}

\TODOAll{There has been some criticism that the ties with enumeration are not strong enough. I disagree. I think we give a solid complicated example, and a Maple template to compute the coefficients. I think that is more powerful than giving static sequences for a handful of cases. If you were interested in strengthening this aspect, I would check the big cactus papers (like the Bergeron that the reviewer was suggesting, Bona and Uhlendick, and see what enumerations they get, localize where it is in OEIS, and give the grammar taht will provide the same result. This would be a way to say "OK, we don't quite have automatic asymptotics yet, but we do rederive everything that has been done bfore in a coherent way." It is a big point, so if you feel pressed for time, feel free to skip it.}
\TODOAll{The last "big" thing to do, would be to tighten the references to past cactus work, and maybe be more explicit about what we mean. (You can try to give it a shot, and I will try to tomorrow afternoon.)}
\TODOAll{Finally, there is a small gimmicky but cute thing you could do, which is the last todonote, which is add a (cleaned up) version of the methodology diagram.}
\bibliographystyle{plain}
\bibliography{bib}{}

% \TODOAll{In Appendix, mention bijection between pure cacti and necklaces.}
\end{document}